\documentclass[11pt]{nyjm}

\title{Generalized improper integral definition for finite limit} 

\author{Michael A. Blischke}  
\address{} 
\email{mblischke@att.net}   



\newtheorem{thm}{Theorem}

\theoremstyle{definition}
\newtheorem{definition}{Definition}
\theoremstyle{example}
\newtheorem{example}{Example}

\usepackage{url}

\begin{document} 
 
\begin{abstract}  
 A generalization of the definition of a one-dimensional improper integral with a finite limit is presented.  The new definition extends the range of valid integrals to include integrals which were previously considered to not be integrable.  This definition is shown to be equivalent to the infinite limit definition presented in ``Generalized improper integral definition for infinite limit'' via a particular change of variable of integration.  The definition preserves linearity and uniqueness.  Integrals which are valid under the conventional definition have the same value under the new definition.  Criteria for interchanging the order of integration and differentiation, and for interchanging the order with a second integration, are obtained.  Examples are provided.
\end{abstract} 

\maketitle
\tableofcontents


\section{Introduction and preliminaries}

There are two basic types of improper integrals. Integrals with an infinite limit are defined as the limit of a series of proper integrals as one of the limits approaches infinity.  Improper integrals with finite limits are needed when the integrand does not have a finite limiting value as the variable of integration approaches a particular ``critical'' value.  In this case, the improper integral is defined as the limit of a series of proper integrals as one of the limits approaches the critical value.  Improper integrals with more than one critical value, or with interior critical values, can be found as a sum of these two basic types.

The improper integral with an infinite upper limit defined by
\[
\int_{a}^{\infty} f(x) \,dx \equiv \lim_{b\rightarrow\infty} \left\{ \int_{a}^{b} f(x) \,dx \right\}
\]
exists when the limit exists, with a similar definition for an integral with an infinite lower limit.

Similarly, the improper integral with a critical lower limit, defined by
\[
\int_{\alpha}^{\beta} g(u) \,du \equiv \lim_{\delta \rightarrow 0^+} \left\{ \int_{\alpha + \delta}^{\beta} g(u) \,du \right\}
\]
exists when that limit exists, with a similar definition for an integral with a critical lower limit.  For any of the above, when the limit does not exist, the integrals are said to not exist, or to diverge.  We will refer to the above as the conventional definitions.

In ``Generalized improper integral definition for infinite limit''~\cite{blischke:ALtDefInfLim2}, the following definition for an integral with an infinite limit was introduced:
\begin{equation}
\label{E001}
\textsf{Z} \hspace{-10.8 bp} \int^{\infty}_{a} f(x) \,dx \equiv \lim_{b\rightarrow\infty}\left\{\int^{b}_{a}f(x) \,dx + \int^{b+c}_{b} f(x) z(x-b) \,dx\right\}
\end{equation}
where $f(x)$ is the function to be integrated, and where $z(x)$ is a termination function, defined therein.  The inclusion of the additional term inside the limit allows convergence to be rigorously shown for a greater range of functions, $f(x)$.  The over-struck \sffamily Z \normalfont on integrals using the alternate definition was included there to distinguish them from integrals that exist using conventional definitions.

A useful form that is equivalent to (\ref{E001}) is
\begin{equation}
\label{E002}
\textsf{Z} \hspace{-10.8 bp} \int^{\infty}_{a} f(x) \,dx = -F(a) - \lim_{b\to\infty}\left\{\int^{c}_{0} F(x+b) z'(x) \,dx\right\}
\end{equation}
where $F(x)$ is defined (for some arbitrary lower limit $\phi$) by
\begin{equation}
\label{E003}
F(x) \equiv \int_{\phi}^{x} f(x') \,dx'.
\end{equation}

When the integrals and the limit in (\ref{E001}) or (\ref{E002}) exist, the integral of $f(x)$ is said to exist under the alternate definition, and to have the value of the limit.
In ~\cite{blischke:ALtDefInfLim2}, it was shown that for all termination functions for which the limit exists, the integral will have the same value, so that the definition gives a unique value.  Many other properties of integrals found using the alternative definition were shown there.

\section{Generalized definition for finite limit improper integral}

It is desirable to introduce a corresponding definition for an improper integral where the critical limit is finite.  We introduce as our definition of improper integral for a function with a critical lower limit $\alpha$ and a noncritical upper limit $\beta$:
\begin{definition}
\begin{equation}
\label{E009}
\textsf{Z} \hspace{-10.8 bp} \int_{\alpha}^{\beta} g(u) \,du \equiv \lim_{\delta \to 0^+} \left\{\int_{\alpha}^{\alpha+\delta}  g(u) w\left(\frac{(u-\alpha)}{\delta}\right) \,du + \int_{\alpha+\delta}^{\beta} g(u) \,du \right\}
\end{equation}
\end{definition}
The term $w(v)$ is defined below.  We will again follow the convention of using the overstruck \sffamily Z \normalfont in this paper for the integrals using the new definitions to distinguish them from conventionally defined integrals.  An improper integral with critical finite upper limit and noncritical lower limit is defined similarly as
\begin{definition}
\begin{equation}
\label{E009b}
\textsf{Z} \hspace{-10.8 bp} \int_{\alpha}^{\beta} g(u) \,du \equiv \lim_{\delta \to 0^+} \left\{\int_{\alpha}^{\beta-\delta} g(u) \,du + \int_{\beta-\delta}^{\beta}  g(u) w\left(\frac{(\beta - u)}{\delta}\right) \,du \right\}.
\end{equation}
\end{definition}

For simplicity in the derivations, and without loss of generality, for the remainder of this paper we will take the critical limit to be the lower limit, and to be 0, giving
\begin{equation}
\label{E010}
\textsf{Z} \hspace{-10.8 bp} \int_{0}^{\beta} g(u) \,du \equiv \lim_{\delta \to 0^+} \left\{\int_{0}^{\delta}  g(u) w(u/\delta) \,du + \int_{\delta}^{\beta} g(u) \,du \right\}.
\end{equation}

The function $w(v)$ will be referred to as the initialization function, analogous to the termination functions in ~\cite{blischke:ALtDefInfLim2} for the infinite limit case.  

As was the case for the termination function, the initialization function is not arbitrary.  Its function is to smooth out the sharp lower bound when taking the limit.  It is required to be finite, to not depend on $\delta$, and to satisfy the following conditions:
\begin{equation}
\label{E011}
w(v) = \left\{ \begin{array}{lcc}
0 \mbox{ } \mbox{ } \mbox{ } v < \epsilon \\
1 \mbox{ } \mbox{ } \mbox{ } v \geq 1 \\
\end{array} \right.
\end{equation}
and
\begin{equation}
\label{E012}
w'(v) \geq 0
\end{equation}
for some $\epsilon$ with
\begin{equation}
\label{E015}
0 < \epsilon \leq 1.
\end{equation}

\begin{definition}
An initialization function is any function satisfying the conditions given in (\ref{E011}) through (\ref{E015}).  
\end{definition}

From (\ref{E011}) we have that
\begin{equation}
\label{E014}
w'(v) = \left\{ \begin{array}{lcc}
0 \mbox{ } \mbox{ } \mbox{ } v < \epsilon \\
0 \mbox{ } \mbox{ } \mbox{ } v > 1 \\
\end{array} \right.
\end{equation}

A less restrictive condition on $w(v)$ may be possible, for example allowing $\epsilon = 0$ with $w(v) \to 0$ suitably fast as $v \to 0$.  However, requiring nonzero $\epsilon$ will allow us to show equivalence between the finite limit definition and the infinite limit definition from ~\cite{blischke:ALtDefInfLim2}.  Investigation into less restrictive conditions on initialization functions which maintain the desirable properties of the new definition, but which expand the allowable set of integrable functions, is a subject for future research.

With (\ref{E011}) and (\ref{E014}), we have that
\begin{equation}
\label{E016}
\int_{0}^{1} w'(v) = \int_{\epsilon}^{1} w'(v) = 1.
\end{equation}
We can get $w(v)$ for $0 < v < 1$ from $w'(v)$ as
\begin{equation}
\label{E017}
w(v) = \int_{0}^{v} w'(v') \,dv'.
\end{equation}

We now define
\begin{equation}
\label{E018}
G(u) \equiv \int_{\phi}^{u} g(s) \,ds
\end{equation}
and using integration by parts, (\ref{E010}) becomes
\begin{eqnarray}
\label{E019}
\textsf{Z} \hspace{-10.8 bp} \int_{0}^{\beta} g(u) \,du \!\! &=& \!\! \lim_{\delta \to 0^+} \left\{ \!\! \left.
\begin{array}{l}
\!\!\!\! \\
\!\!\!\! \end{array}
G(u) w(u/\delta) \right|_{0}^{\delta} - \frac{1}{\delta} \int_{0}^{\delta}  G(u) w'(u/\delta) \,du + \!\! \left.
\begin{array}{l}
\!\!\!\! \\
\!\!\!\! \end{array}
G(u) \right|_{\delta}^{\beta} \right\} \nonumber \\
\!\! &=& \!\! G(\beta) -  \lim_{\delta \to 0^+} \left\{ \frac{1}{\delta} \int_{0}^{\delta}  G(u) w'(u/\delta) \,du \right\}.
\end{eqnarray}

This form is analogous to (\ref{E002}), and can be an easier form to work with than (\ref{E010}).

We will follow a path analogous to that taken with termination functions in ~\cite{blischke:ALtDefInfLim2}.  Given two initialization functions $w_1(v)$ and $w_2(v)$, we can combine them to obtain a third, via their derivatives, as
\begin{equation}
\label{E020}
w'(v) \equiv \int_0^1 \frac{w_1 ' \left( v / v' \right) w_2 ' (v')}{v'} \,dv'.
\end{equation}
Using (\ref{E014}), (\ref{E020}) can also be written as
\begin{equation}
\label{E021}
w'(v) = \int_v^1 \frac{w_1 ' \left( v / v' \right) w_2 ' (v')}{v'} \,dv'.
\end{equation}

From (\ref{E020}), since $w_1 ' (v)$ and $w_2 ' (v)$ are both non-negative, $w ' (v)$ is also non-negative, satisfying (\ref{E012}).  If $w_1 ' (v)$ satisfies (\ref{E014}), we can see that $w'(v)$ also satisfies (\ref{E014}).  We can also show that $w'(v)$ satisfies (\ref{E016}), and therefore also (\ref{E011}):

\begin{eqnarray}
\label{E022}
\int_0^1 w'(v) \!\! &=& \!\! \int_0^1 \left[ \int_0^1 \frac{w_1 ' \left( v / v' \right) w_2 ' (v')}{v'} \,dv' \right]  \,dv \nonumber \\
\!\! &=& \!\! \int_0^1 \left[ \int_0^1 w_1 ' \left( v / v' \right) \,dv \right] \frac{ w_2 ' (v')}{v'} \,dv' \nonumber \\
\!\! &=& \!\! \int_0^1 \left[ \int_0^{v'} w_1 ' \left( v / v' \right) \,dv \right] \frac{ w_2 ' (v')}{v'} \,dv' \nonumber \\
\!\! &=& \!\! \int_0^1 v' \frac{ w_2 ' (v')}{v'} \,dv' \nonumber \\
\!\! &=& \!\! 1.
\end{eqnarray}

We will denote the combined initialization function using the same notation we used for termination functions,
\begin{equation}
\label{E023}
w(v) = w_1(v) \odot w_2(v).
\end{equation}

We also have, substituting $v'' = v / v'$, that
\begin{eqnarray}
\label{E024}
w'(v) = \int_0^1 \frac{w_1 ' \left( v / v' \right) w_2 ' (v')}{v'} \,dv' \!\! &=& \!\! \int_0^1 \frac{w_1 ' \left( v'' \right) w_2 ' (v/v'')}{v'} \frac{v}{v' v''^2} \,dv' \nonumber \\
\!\! &=& \!\! \int_0^1 \frac{w_1 ' \left( v'' \right) w_2 ' (v/v'')}{v''} \,dv'
\end{eqnarray}
so the relation (\ref{E023}) satisfies commutivity.

\section{Conversion to infinite limit integral}

We will now examine what happens when we switch between the generalized finite limit integral and an infinite limit improper integral. This development sets up the next section, where we show equivalence between the generalized finite limit definition and the generalized infinite limit definition.  

We will use the general change of variable defined by
\begin{eqnarray}
\label{E025}
u \!\! &=& \!\! \psi(x) \\
x \!\! &=& \!\! \psi^{-1}(u) . \nonumber
\end{eqnarray}
We require for all $x > \psi^{-1}(\beta)$ that $\psi(x)$ be finite, that it be strictly monotonic with
\begin{equation}
\label{E026}
\psi'(x) < 0,
\end{equation}
and that
\begin{equation}
\label{E027}
\lim_{x \to \infty} \left\{ \psi(x) \right\} \equiv 0.
\end{equation}
Combining these properties, we obtain that
\begin{equation}
\label{E027b}
\psi(x) > 0,
\end{equation}
that
\begin{equation}
\label{E028}
0 \leq \frac{\psi(x+b)}{\psi(b)} \leq 1 \mbox{  } \mbox{ for } \mbox{  }  b \geq \psi^{-1}(\beta) \mbox{ , } x \geq 0
\end{equation}
and that
\begin{equation}
\label{E02bb}
\int_{x}^{\infty} \psi'(x') dx' = \psi(x).
\end{equation}

Using (\ref{E025}) in (\ref{E018}) gives
\begin{equation}
\label{E029}
G(u) = \int_{\psi^{-1}(\phi)}^{\psi^{-1}(u)} g(\psi(t)) \psi'(t) \,dt.
\end{equation}
Defining
\begin{equation}
\label{E030}
f(t) \equiv -g(\psi(t)) \psi'(t)
\end{equation}
and
\begin{equation}
\label{E031}
F(x) \equiv \int_{\psi^{-1}(\phi)}^{x} f(t) \,dt
\end{equation}
we get
\begin{equation}
\label{E032}
G(u) = -\int_{\psi^{-1}(\phi)}^{\psi^{-1}(u)} f(t) \,dt = -F(\psi^{-1}(u))
\end{equation}
or
\begin{equation}
\label{E033}
F(x) = -G(\psi(x)).
\end{equation}

We will now define an infinite limit integral corresponding to our finite limit integral as
\begin{equation}
\label{E034}
\Xi \hspace{-11.2 bp} \int_{\psi^{-1}(\beta)}^{\infty} f(x) \,dx \equiv \textsf{Z} \hspace{-10.8 bp} \int_{0}^{\beta} g(u) \,du.
\end{equation}
In (\ref{E034}) we use an over-struck $\Xi$, instead of an over-struck \sffamily Z\normalfont, to allow us to distinguish the resulting integral (\ref{E042}) from the infinite limit integral definition introduced in ~\cite{blischke:ALtDefInfLim2}.  We will use the overstruck $\Xi$ throughout this paper for the infinite limit improper integral that is obtained from our finite limit improper integral via a change of variable \emph{within the definition}.  Since the definition is in terms of conventional integrals, we know we can safely perform those changes of variable.  Note that we haven't yet shown any relation beween this and the infinite limit integral definition from ~\cite{blischke:ALtDefInfLim2}.

Using (\ref{E034}) in (\ref{E019}) we get
\begin{eqnarray}
\label{E035}
\Xi \hspace{-11.2 bp} \int_{\psi^{-1}(\beta)}^{\infty} f(x) \,dx \!\! &=& \!\! G(\beta) -  \lim_{\delta \to 0^+} \left\{ \frac{1}{\delta} \int_{0}^{\delta}  G(u) w'(u/\delta) \,du \right\} \nonumber \\
&=& \!\! G(\beta) -  \lim_{\delta \to 0^+} \left\{ \frac{1}{\delta} \int_{\psi^{-1}(\delta)}^{\infty} G(\psi(x)) w'\left(\frac{\psi(x)}{\delta}\right) \psi'(x) \,dx \right\} \nonumber \\
&=& \!\! -F\left(\psi^{-1}(\beta)\right) + \lim_{\delta \to 0^+} \left\{ \frac{1}{\delta} \int_{\psi^{-1}(\delta)}^{\infty}  F(x) w'\left(\frac{\psi(x)}{\delta}\right) \psi'(x) \,dx \right\} \nonumber \\
&=& \!\! -F\left(\psi^{-1}(\beta)\right) + \lim_{b \to \infty} \left\{ \frac{1}{\psi(b)} \int_{b}^{\infty} F(x) w'\left(\frac{\psi(x)}{\psi(b)}\right) \psi'(x) \,dx \right\} \nonumber \\
& &
\end{eqnarray}
where $b \equiv \psi^{-1}(\delta)$.

We now define $\zeta(x,b)$ such that
\begin{equation}
\label{E036}
\zeta(0,b) \equiv 1
\end{equation}
and
\begin{equation}
\label{E037}
\zeta'(x-b,b) \equiv \left\{ \begin{array}{cc}
{w'\left(\frac{\psi(x)}{\psi(b)}\right) \frac{\psi'(x)}{\psi(b)}} & x \geq b \\
0 & x < b \\
\end{array} \right.
\end{equation}
or
\begin{equation}
\label{E038}
\zeta'(x,b) = \left\{ \begin{array}{cc}
{w'\left(\frac{\psi(x+b)}{\psi(b)}\right) \frac{\psi'(x+b)}{\psi(b)}} & x \geq 0 \\
0 & x < 0 \\
\end{array} \right. .
\end{equation}
We see from (\ref{E012}), (\ref{E026}), and (\ref{E027b}) that 
\begin{equation}
\label{E038b}
\zeta'(x,b) \leq 0
\end{equation}
for all $x$.

From (\ref{E014}), $\zeta'(x,b) = 0$ when $\psi(x+b) < \epsilon \psi(b)$.  We will define $e(\epsilon,b)$ by
\begin{equation}
\label{E038c}
\frac{\psi(e(\epsilon,b)+b)}{\psi(b)} \equiv \epsilon
\end{equation}
so we get
\begin{equation}
\label{E038d}
e(\epsilon,b) = \psi^{-1}\left( \epsilon \psi(b)\right) - b
\end{equation}
and then
\begin{equation}
\label{E041}
\zeta'(x,b) = 0 \mbox{ } \mbox{ for } x > e(\epsilon,b) .
\end{equation}

Since $w'(v)$ and $\frac{\psi'(x+b)}{\psi(b)}$ are finite, we have that $\zeta'(x,b)$ is finite, and from (\ref{E038}) and (\ref{E041}) it is seen that
\begin{equation}
\label{E040}
\int_{-\infty}^{\infty}  \zeta'(x,b) \,dx = \int_{0}^{e(\epsilon,b)}  \zeta'(x,b) \,dx = -1.
\end{equation}
The integral in (\ref{E040}) is proper and absolutely convergent.  

In a reciprocal fashion, if we have $\zeta(x,b)$ that can be shown to satisfy (\ref{E036}) through (\ref{E038b}) and (\ref{E041}) for some $w'(v)$, with $\zeta'(x,b)$ finite, then (\ref{E011}) through (\ref{E016}) can be satisfied.

Substituting, (\ref{E035}) becomes
\begin{eqnarray}
\label{E042}
\Xi \hspace{-11.2 bp} \int_{\psi^{-1}(\beta)}^{\infty} f(x) \,dx \!\! &=& \!\! -F\left(\psi^{-1}(\beta)\right) + \lim_{b \to \infty} \left\{ \int_{b}^{\infty} F(x) \zeta'(x-b,b) \,dx \right\} \nonumber \\
&=& \!\! -F\left(\psi^{-1}(\beta)\right) + \lim_{b \to \infty} \left\{  \int_{0}^{\infty} F(x+b) \zeta'(x,b) \,dx \right\}.
\end{eqnarray}

Because of the similarity of (\ref{E042}) and (\ref{E002}), we will refer to $\zeta (x,b)$ also as a termination function, even though it is not constant WRT $b$.  When it is not obvious from context which termination function we are referring to, we will call the $\zeta (x,b)$ a termination function of the second type, and will call termination functions as described in ~\cite{blischke:ALtDefInfLim2} a termination function of the first type.

Using the relations between $w(v)$ and $\zeta(x,b)$, we can freely switch between the finite limit integral (\ref{E019}) and its corresponding infinite limit integral (\ref{E042}).  We'll denote the relation between an initialization function $w(v)$ used in the finite limit integral and the termination function $\zeta(x,b)$ used in the corresponding infinite limit integral and given by (\ref{E037}) or (\ref{E038}), as
\begin{equation}
\label{E043}
\zeta(x) \Leftrightarrow w(v).
\end{equation}

We next show that combining initialization functions and combining termination functions of the second type are equivalent.  Beginning with (\ref{E021}) and using
\begin{equation}
\label{E044}
v = \frac{\psi(x+b)}{\psi(b)}
\end{equation}
\begin{equation}
\label{E045}
v' = \frac{\psi(x'+b)}{\psi(b)}
\end{equation}
we get
\begin{equation}
\label{E046}
dv' = \frac{\psi'(x'+b)}{\psi(b)} \,dx'
\end{equation}
and
\begin{eqnarray}
\label{E047}
w'\left(\frac{\psi(x+b)}{\psi(b)}\right) \!\! &=& \!\! -\int_0^x \frac{w_1 ' \left(\frac{\psi(x+b)}{\psi(x'+b)}\right) w_2 ' \left(\frac{\psi(x'+b)}{\psi(b)}\right)}{\left(\psi(x'+b) / \psi(b)\right)} \frac{\psi'(x'+b)}{\psi(b)} \,dx' \nonumber \\
\!\! &=& \!\! -\int_0^x {\frac{\left[ w_1 ' \left(\frac{\psi(x+b)}{\psi(x'+b)}\right) \frac{\psi'(x'+b)}{\psi(x'+b)} \right] \left[ w_2 '\left(\frac{\psi(x'+b)}{\psi(b)}\right) \frac{\psi'(x'+b)}{\psi(b)} \right]}{\left(\psi'(x+b) / \psi(b)\right)} } \,dx' \nonumber \\
\!\! &=& \!\! -\int_b^{x+b} {\frac{\left[ w_1 ' \left(\frac{\psi(x+b)}{\psi(x')}\right) \frac{\psi'(x')}{\psi(x')} \right] \left[ w_2 ' \left(\frac{\psi(x')}{\psi(b)}\right) \frac{\psi'(x')}{\psi(b)} \right]}{\left(\psi'(x+b) / \psi(b)\right)}}  \,dx' .
\end{eqnarray}
Using (\ref{E037}) and (\ref{E038})
\begin{eqnarray}
\label{E048}
w'\left(\frac{\psi(x+b)}{\psi(b)}\right) \!\! &=& \!\! -\frac{\psi(b)}{\psi'(x+b)} \int_b^{x+b} \zeta'_1(x+b-x',b) \zeta'_2(x'-b,b) \,dx' \nonumber \\
\!\! &=& \!\! -\frac{\psi(b)}{\psi'(x+b)} \int_0^{x} \zeta'_1(x-x',b) \zeta'_2(x',b) \,dx'
\end{eqnarray}
and so
\begin{eqnarray}
\label{E049}
\zeta'(x,b) = -\int_0^x \zeta'_1(x-x',b) \zeta'_2(x',b) \,dx'
\end{eqnarray}

Thus we have that the combination of two termination functions of the second type is also a termination function of the second type.  Note that each correspondance between termination function and initialization function uses the \emph{same} change of variable, $u=\psi(x)$.

We write this more compactly as
\begin{equation}
\label{E050}
\zeta'(x,b) = \zeta'_1(x,b) \otimes \zeta'_2(x,b)
\end{equation}
and denote the combined termination function using the notation
\begin{equation}
\label{E051}
\zeta(x,b) = \zeta_1(x,b) \odot \zeta_2(x,b).
\end{equation}

Thus, given
\begin{equation}
\label{E052}
\zeta_1(x,b) \Leftrightarrow w_1(v) \mbox{  } \mbox{ and } \mbox{  } \zeta_2(x,b) \Leftrightarrow w_2(v)
\end{equation}
we find that
\begin{equation}
\label{E053}
\zeta_1(x,b) \odot \zeta_2(x,b) \Leftrightarrow w_1(v) \odot w_2(v).
\end{equation}

That is, given two initialization functions, the combination of their corresponding termination functions is equal to the termination function corresponding to their combination.

\section{Equivalence of generalized definitions}

We are now ready to show equivalence between our finite limit and infinite limit definitions.  Equation (\ref{E042}) is almost the same form as (\ref{E002}) and Eq. (9) of ~\cite{blischke:ALtDefInfLim2}, the only difference being that $\zeta(x,b)$ is a function of $b$, unlike the termination functions described in ~\cite{blischke:ALtDefInfLim2}.  For a particular choice of $\psi(x)$, however, the dependence of $\zeta(x,b)$ on $b$ vanishes.  In that case, the infinite limit integral corresponding to our finite limit improper integral, (\ref{E042}), is identical to the infinite limit integral presented in ~\cite{blischke:ALtDefInfLim2}.

Choosing, with $\alpha > 0$,
\begin{equation}
\label{E054}
\psi(x) = e^{-\alpha x}
\end{equation}
we get
\begin{equation}
\label{E055}
\psi'(x) = -\alpha e^{-\alpha x}.
\end{equation}
Substituting these into (\ref{E038}) we can write
\begin{equation}
\label{E056}
\zeta'(x,b) = -\alpha w'\left(e^{-\alpha x}\right) e^{-\alpha x} \equiv z'(x).
\end{equation}
Since this is not a function of $b$, it satisfies the requirements for a termination function.  Using this transformation, we can see that for every finite limit improper integral using the general definition (\ref{E010}), there is a corresponding infinite-limit improper integral using the definition from ~\cite{blischke:ALtDefInfLim2}.  We also have that for any termination function $z(x)$, the initialization function can be explicitly found as
\begin{equation}
\label{E057}
w'\left(e^{-\alpha x}\right) = -z'(x) \frac{e^{\alpha x}}{\alpha}
\end{equation}
\begin{equation}
\label{E058}
w'\left(u\right) =  \frac{-z'\left(\frac{-ln(u)}{\alpha}\right)}{\alpha u}.
\end{equation}
We thus have, when the change of variable (\ref{E025}) is given by (\ref{E054}), that
\begin{equation}
\label{E059}
\Xi \hspace{-11.2 bp} \int_{\frac{-ln(\beta)}{\alpha}}^{\infty} f(x) \,dx = \textsf{Z} \hspace{-10.8 bp} \int_{\frac{-ln(\beta)}{\alpha}}^{\infty} f(x) \,dx
\end{equation}
and the definition for the infinite limit case is seen to be equivalent to the definition for the finite limit case, so
\begin{equation}
\label{E060}
\textsf{Z} \hspace{-10.8 bp} \int_{0}^{\beta} g(u) \,du = \textsf{Z} \hspace{-10.8 bp} \int_{\frac{-ln(\beta)}{\alpha}}^{\infty} f(x) \,dx
\end{equation}

We find that $c$ and $\epsilon$ are related as
\begin{equation}
\label{E061}
\epsilon = e^{-\alpha c}.
\end{equation}

The equivalence between the finite limit and infinite limit generalized definitions given by (\ref{E054}) and the correspondence between initialization and termination functions given by (\ref{E055}) means that the properties found for the infinite limit case in ~\cite{blischke:ALtDefInfLim2} all have corresponding properties in the finite limit case.  We will list the properties here.

If $w_1(\nu)$ is an initialization function for $g(u)$, then for any other initialization function $w_2(\nu)$, $w(\nu)$ given by (\ref{E020}) is also an initialization function for $g(u)$, and gives the same value for the integral.  

The integral defined using (\ref{E010}) produces a unique value for all initialization functions for which the limit exists.  

When the integral exists using the conventional definition, our general definition gives the same answer,

\begin{equation}
\label{E063}
\textsf{Z} \hspace{-10.8 bp} \int_{0}^{\beta} g(u) \,du  = \int_{0}^{\beta} g(u) \,du.
\end{equation}

Our general definition satisfies linearity,
\begin{equation}
\label{E064}
a \textsf{Z} \hspace{-10.8 bp} \int_{0}^{\beta} g(u) \,du + b \textsf{Z} \hspace{-10.8 bp} \int_{0}^{\beta} h(u) \,du = \textsf{Z} \hspace{-10.8 bp} \int_{0}^{\beta} \left[ a g(u) + b h(u) \right] \,du
\end{equation}
when the integrals on the LHS both exist.

Differentiation under the integral sign can be performed when both of the integrals exist, and we have
\begin{equation}
\label{E065}
\frac{d}{dy} \left[ \textsf{Z} \hspace{-10.8 bp} \int_{0}^{\beta} g(u,y) \,du \right] = \left[ \textsf{Z} \hspace{-10.8 bp} \int_{0}^{\beta} \frac{\partial}{\partial y} g(u,y) \,du \right].
\end{equation}
When one or the other integrals exists with the conventional definition, we also get
\begin{equation}
\label{E066}
\frac{d}{dy} \left[ \int_{0}^{\beta} g(u,y) \,du \right] = \left[ \textsf{Z} \hspace{-10.8 bp} \int_{0}^{\beta} \frac{\partial}{\partial y} g(u,y) \,du \right]
\end{equation}
and
\begin{equation}
\label{E067a}
\frac{d}{dy} \left[ \textsf{Z} \hspace{-10.8 bp} \int_{0}^{\beta} g(u,y) \,du \right] = \left[ \int_{0}^{\beta} \frac{\partial}{\partial y} g(u,y) \,du \right].
\end{equation}

Interchange of the order of iterated integrations is also allowed.  Here we assume that
\begin{equation}
\label{E067b}
\textsf{Z} \hspace{-10.8 bp} \int_{0}^{\beta} g(u,y) \,dx
\end{equation}
exists over the domain $\gamma \leq y \leq \delta$ for some initialization function $w(v,y)$, and that
\begin{equation}
\label{E068a}
h(u,t) \equiv \int^{t} g(u,y) s(y) \,dy
\end{equation}
exists, using the Riemann definition, over the domain $\gamma \leq t \leq \delta$, for $0 < u \leq \beta$.  The function $s(y)$ is arbitrary.  We further assume that
\begin{equation}
\label{E068b}
\textsf{Z} \hspace{-10.8 bp} \int_{0}^{\beta} h(u,y) \,dx
\end{equation}
exists over the domain $\gamma \leq y \leq \delta$ for some initialization function $\widetilde{w}(v,y)$.

\begin{thm}
When integrals in (\ref{E067b}) and (\ref{E068b}) exist, and the Riemann integral in (\ref{E068a}) exists,
 \begin{equation}
\label{E071}
\int_{\gamma}^{\delta} s(y) \left[\textsf{Z} \hspace{-10.8 bp} \int_{0}^{\beta} g(u,y) \,du \right] dy = \textsf{Z} \hspace{-10.8 bp} \int_{0}^{\beta} \left[ \int_{\gamma}^{\delta} s(y) g(u,y) \,dy \right] \,du.
\end{equation}
\end{thm}

\begin{proof}
From ~\cite{blischke:ALtDefInfLim2} and (\ref{E059}), we have that
\begin{equation}
\label{E069}
\int_{\gamma}^{\delta} s(y) \left[\Xi \hspace{-11.2 bp} \int_{\frac{-ln(\beta)}{\alpha}}^{\infty} f(x,y) \,dx \right] dy = \Xi \hspace{-11.2 bp} \int_{\frac{-ln(\beta)}{\alpha}}^{\infty} \left[ \int_{\gamma}^{\delta} s(y) f(x,y) \,dy \right] \,dx.
\end{equation}

Using (\ref{E068a}) and (\ref{E059}), we get
\begin{eqnarray}
\label{E070}
\Xi \hspace{-11.2 bp} \int_{\frac{-ln(\beta)}{\alpha}}^{\infty} \left[ \int_{\gamma}^{\delta} s(y) f(x,y) \,dy \right] \,dx \!\! &=& \!\! \Xi \hspace{-11.2 bp} \int_{\frac{-ln(\beta)}{\alpha}}^{\infty} \left[ \int_{\gamma}^{\delta} -s(y) g(\psi(x),y)\psi'(x) \,dy \right] \,dx \nonumber \\
\!\! &=& \!\! \Xi \hspace{-11.2 bp} \int_{\frac{-ln(\beta)}{\alpha}}^{\infty} \left. \stackrel{}-\!\!h(\psi(x),y)\psi'(x) \right|_{\gamma}^{\delta} \,dx \nonumber \\
\!\! &=& \!\! \textsf{Z} \hspace{-10.8 bp} \int_{0}^{\beta} \left. \stackrel{}{h(u,y)} \right|_{\gamma}^{\delta} \,du \nonumber \\
\!\! &=& \!\! \textsf{Z} \hspace{-10.8 bp} \int_{0}^{\beta} \left[ \int_{\gamma}^{\delta} s(y) g(u,y) \,dy \right] \,du.
\end{eqnarray}
Using (\ref{E059}) and (\ref{E060}),
\begin{equation}
\label{E071b}
\int_{\gamma}^{\delta} s(y) \left[\Xi \hspace{-11.2 bp} \int_{\frac{-ln(\beta)}{\alpha}}^{\infty} f(x,y) \,dx \right] dy = \int_{\gamma}^{\delta} s(y) \left[\textsf{Z} \hspace{-10.8 bp} \int_{0}^{\beta} g(u,y) \,du \right] dy.
\end{equation}
\end{proof}

When one side exists using the conventional definition, we also get either
\begin{equation}
\label{E072}
\int_{\gamma}^{\delta} s(y) \left[\textsf{Z} \hspace{-10.8 bp} \int_{0}^{\beta} g(u,y) \,du \right] dy = \int_{0}^{\beta} \left[ \int_{\gamma}^{\delta} s(y) g(u,y) \,dy \right] \,du
\end{equation}
or
\begin{equation}
\label{E073}
\int_{\gamma}^{\delta} s(y) \left[ \int_{0}^{\beta} g(u,y) \,du \right] dy = \textsf{Z} \hspace{-10.8 bp} \int_{0}^{\beta} \left[ \int_{\gamma}^{\delta} s(y) g(u,y) \,dy \right] \,du.
\end{equation}

A change of variable of integration of the form $u' = c u$ for nonzero constant $c$ is valid.  An arbitrary change of the variable of integration is not necessarily valid.

\section{Examples}

The following two examples show the evaluation of integrals which do not exist using the conventional definition.

\begin{example}
\label{EX01}
\[
g(u) = \frac{\sin(1/u)}{u^2}
\]
\[
G(u) = \int_{0}^{u} \frac{\sin(1/s)}{s^2} \,ds = \cos(1/u)
\]

Using (\ref{E019}) we have
\[
\textsf{Z} \hspace{-10.4 bp} \int_{0}^{1/a} \frac{\sin(1/u)}{u^2} \,du = G(1/a) -  \lim_{\delta \to 0^+} \left\{ \frac{1}{\delta} \int_{0}^{\delta}  G(u) w'(u/\delta) \,du \right\}.
\]

We can use
\[
w(v) = \left\{ \begin{array}{ccc}
{(2v-1)} \mbox{ } \mbox{ } \mbox{ } {1/2 \leq v \leq 1} \\
{\mbox{ } \mbox{ } \mbox{ } 0} {\mbox{ } \mbox{ } \mbox{ } \mbox{ } \mbox{ } \mbox{ } \mbox{ } \mbox{ } \mbox{ } \mbox{ } \mbox{ } \mbox{ } \mbox{ } \mbox{ }} {v < 1/2} \\
\end{array} \right.
\]
so
\[
w'(v) = \left\{ \begin{array}{ccc}
{2} {\mbox{ } \mbox{ } \mbox{ } \mbox{ } \mbox{ } \mbox{ }} {1/2 \leq v \leq 1} \\
{0} {\mbox{ } \mbox{ } \mbox{ } \mbox{ } \mbox{ } \mbox{ } \mbox{ } \mbox{ } \mbox{ } \mbox{ } \mbox{ } \mbox{ } \mbox{ } \mbox{ }} \mbox{ else } \\
\end{array} \right.
\]
giving us
\[
\textsf{Z} \hspace{-10.4 bp} \int_{0}^{1/a} \frac{\sin(1/u)}{u^2} \,du = \cos(a) -  \lim_{\delta \to 0^+} \left\{ \frac{1}{\delta} \int_{\delta/2}^{\delta}  2 \cos(1/u) \,du \right\} . \nonumber
\]

Evaluating the integral on the RHS gives
\[
\int_{\delta/2}^{\delta}  \cos(1/u) \,du = \left. \frac{-2}{\delta}\left( u \cos(1/u) + \mbox{Si}(1/u) \right) \right|_{\delta/2}^{\delta} \nonumber
\]
where $Si$ is the sine integral.  Expanding this to the necessary order of argument in the large argument ($u \to 0^+$) limit,
\[
\mbox{Si}(1/u) \approx \frac{\pi}{2} - u\cos(1/u)-u^2\sin(1/u).
\]
The cosine terms cancel, and the $\pi/2$ term doesn't contribute when the limits are taken.  Thus, the integral is of order $\delta^2$, and we get for the limit
\[
\lim_{\delta \to 0^+} \left\{ O(\delta) \right\} = 0
\]
so we obtain that
\[
\textsf{Z} \hspace{-10.4 bp} \int_{0}^{1/a} \frac{\sin(1/u)}{u^2} \,du = \cos(a).
\]
\end{example}
\begin{example}
\label{EX02}
\[
g(u) = \frac{\cos(1/u)}{u^3}
\]
\[
G(u) = \int_{0}^{u} \frac{\cos(1/s)}{s^3} \,ds = -\cos(1/u) - \frac{\sin(1/u)}{u}
\]

Using (\ref{E019}) we have
\[
\textsf{Z} \hspace{-10.4 bp} \int_{0}^{1/a} \frac{\cos(1/u)}{u^3} \,du = G(1/a) -  \lim_{\delta \to 0^+} \left\{ \frac{1}{\delta} \int_{0}^{\delta}  G(u) w'(u/\delta) \,du \right\}.
\]

We will use
\[
w(v) = \left\{ \begin{array}{ccc}
{3(2v-1)^2 - 2(2v-1)^3} \mbox{ } \mbox{ } \mbox{ } {1/2 \leq v \leq 1} \\
0 {\mbox{ } \mbox{ } \mbox{ } \mbox{ } \mbox{ } \mbox{ } \mbox{ } \mbox{ } \mbox{ } \mbox{ } \mbox{ } \mbox{ } \mbox{ } \mbox{ } \mbox{ } \mbox{ } \mbox{ } \mbox{ } \mbox{ } \mbox{ } \mbox{ } \mbox{ } \mbox{ } \mbox{ } \mbox{ } \mbox{ } \mbox{ }} {v < 1/2} \\
\end{array} \right.
\]
so
\[
w'(v) = \left\{ \begin{array}{ccc}
{12(2v-1)(2-2v)} {\mbox{ } \mbox{ } \mbox{ }} {1/2 \leq v \leq 1} \\
{0} {\mbox{ } \mbox{ } \mbox{ } \mbox{ } \mbox{ } \mbox{ } \mbox{ } \mbox{ } \mbox{ } \mbox{ } \mbox{ } \mbox{ } \mbox{ } \mbox{ } \mbox{ } \mbox{ } \mbox{ } \mbox{ } \mbox{ } \mbox{ } \mbox{ }} \mbox{ else } \\
\end{array} \right.
\]
giving
\begin{eqnarray}
& & \!\!\!\!\!\! \textsf{Z} \hspace{-10.4 bp} \int_{0}^{1/a} \frac{\cos(1/u)}{u^3} \,du = -\cos(a) - a \sin(a) + \nonumber \\
& & \mbox{   } \lim_{\delta \to 0^+} \left\{ \frac{12}{\delta} \int_{\delta/2}^{\delta}  \left(\cos(1/u) + \frac{\sin(1/u)}{u}\right) (2u/\delta-1)(2-2u/\delta) \,du \right\} . \nonumber
\end{eqnarray}

Evaluating the integral on the RHS gives
\begin{eqnarray}
 & & \!\!\!\!\!\! \int_{\delta/2}^{\delta}  \left(\cos(1/u) + \frac{\sin(1/u)}{u}\right) (2u/\delta-1)(2-2u/\delta) \,du \mbox{   } \mbox{   } \mbox{   } \mbox{   } \mbox{   } \mbox{   }\nonumber \\
& & \mbox{   } \mbox{   } \mbox{   } = \frac{-2}{\delta^2}\left[ \frac{u}{6}\left(6 \delta^2 - 9 u \delta + 4 u^2 + 4 \right) \cos(1/u) \right. + \nonumber \\
& & \mbox{   } \mbox{   } \mbox{   } \mbox{   } \mbox{   } \mbox{   } \mbox{   } \mbox{   } \mbox{   } \mbox{   } \mbox{   } \mbox{   } \mbox{   } \mbox{   } \!\! \left. \frac{u}{6} \left( 4u - 9 \delta \right) \sin(1/u) + \frac{3\delta}{2} \mbox{Ci}(1/u) + \frac{2}{3} \mbox{Si}(1/u) \right]_{\delta /2}^{\delta} \nonumber
\end{eqnarray}
where $Si$ and $Ci$ are the sine integral and cosine integral.  Expanding these to the necessary order in the argument in the large argument ($u \to 0^+$) limit,
\[
\mbox{Ci}(1/u) \approx u\sin(1/u)-u^2\cos(1/u)-2u^3\sin(1/u)
\]
\[
\mbox{Si}(1/u) \approx \frac{\pi}{2} - u\cos(1/u)-u^2\sin(1/u)+2u^3\cos(1/u).
\]
With a little algebra there is much cancellation of terms, and the integral can be shown to be at least of order $\delta^2$.  We thus get for the limit
\[
\lim_{\delta \to 0^+} \left\{ O(\delta) \right\} = 0
\]
so we get
\[
\textsf{Z} \hspace{-10.4 bp} \int_{0}^{1/a} \frac{\cos(1/u)}{u^3} \,du = -\cos(a) - a \sin(a).
\]
 \end{example}

\section{Conclusion}

A generalized definition for an improper integral with finite bounds has been presented.  The definition presented here is a more powerful alternative to the conventional definition.  The range of functions which are integrable under this definition is expanded as compared with the conventional definition.

The new definition presented here gives the same result as the conventional definition when that applies, and preserves uniqueness and linearity.  The generalized definition allows interchange of the order of differentiation and integration whenever the two integrals exist under the new definition.  Also allowed is interchange of the order of integration of iterated integration, again when the integrations exist under this definition.  The ability to rigorously interchange order of integrations, or order of integration and differentiation, in cases where integrals under the conventional definition do not converge, provides an added tool for manipulation of complicated integrals.  An arbitrary change of the variable of integration is not necessarily valid, although scaling the variable of integration by a constant is.

The generalized finite limit definition presented here has been shown to be equivalent to the generalized infinite limit definition presented in ~\cite{blischke:ALtDefInfLim2}.  For a particular change of variable transforming between finite limit and infinite limit integrals, the existence of either the finite or the infinite limit integral implies existence of the other, with the same value.

\bibliographystyle{alpha}
\bibliography{ExtendedIntegration}

\end{document}